\documentclass{svproc}

\usepackage{url}

\usepackage{amssymb}
\usepackage{amsmath}
\usepackage{amsthm}

\usepackage{todonotes}

\begin{document}
\mainmatter
\title{Quaternary Legendre Pairs}
\titlerunning{Quaternary Legendre Pairs}

\author{Ilias S. Kotsireas\inst{1} \and Arne Winterhof\inst{2}}
\authorrunning{Ilias S. Kotsireas and Arne Winterhof}
\tocauthor{Ilias S. Kotsireas, Arne Winterhof}

\institute{CARGO Lab, Wilfrid Laurier University, Waterloo, ON N2L 3C5,
Canada,\\
\email{ikotsire@wlu.ca}\and
RICAM, Austrian Academy of Sciences, Altenbergerstr.\ 69, 4040~Linz, Austria,\\
\email{arne.winterhof@oeaw.ac.at}}


\maketitle

\begin{abstract}
 We introduce quaternary Legendre pairs of length $\ell$. In contrast to binary Legendre pairs they can exist for even $\ell$ as well. 
 First we show that they are pertinent to the construction of quaternary Hadamard matrices of order $2\ell+2$
 and thus of binary Hadamard matrices of order $4\ell+4$.
  Then for a prime $p>2$ we 
 present a construction of a pair of sequences of length $p$ from which we can derive quaternary Legendre pairs of length~$\ell=2p$ by decompression for $p=3,5,7,13,19,31,41$. Moreover, we give also constructions of Legendre pairs of length $\ell$ for all remaining even~$\ell\le 24$.
\end{abstract}


\section{Introduction}

Let $$A=[a_0,a_1,\ldots,a_{\ell-1}]\in \mathbb{C}^\ell$$ be a finite complex sequence of length $\ell$ periodically continued with period $\ell$, that is, 
$$a_{j+\ell}=a_j,\quad j=0,1,\ldots$$
The {\em periodic autocorrelation function} $PAF(A,s)$ of $A$ at lag $s$ is
$$PAF(A,s)=\sum_{j=0}^{\ell-1} a_j\overline{a_{j+s}},\quad s=0,1,\ldots,\ell-1,$$
where $\overline{a}$ is the complex conjugate of $a\in \mathbb{C}$.
Two sequences 
$$A=[a_0,a_1,\ldots,a_{\ell-1}]\quad \mbox{and}\quad B=[b_0,b_1,\ldots,b_{\ell-1}]$$ of length~$\ell$ 
form a {\em Legendre pair} $(A,B)$ if
\begin{equation}\label{Legpairdef}PAF(A,s)+PAF(B,s)=-2,\quad s=1,2,\ldots,\left\lfloor\frac{\ell}{2}\right\rfloor.
\end{equation}
Note that 
$$PAF(A,\ell-s)=\overline{PAF(A,s)},\quad s=1,2,\ldots,\left\lfloor \frac{\ell}{2}\right\rfloor,$$ 
and we may extend the range of $s$ in \eqref{Legpairdef} to $s=1,2,\ldots,\ell-1$.

A Legendre pair $(A,B)$ is {\em binary} if both sequences $A$ and $B$ are binary, that is, $A,B\in \{-1,1\}^\ell$.
Binary Legendre pairs of length $\ell$ are pertinent to the construction of Hadamard matrices of order $2\ell+2$, see for example \cite[Theorem 3]{FGS2001}, \cite[Theorem 1]{K2013} and \cite[Theorem~2]{ABH2020}. Since the order of an Hadamard matrix is either $2$ or divisible by $4$, there is no binary Legendre pair $(A,B)$ of even length~$\ell$.
Binary Legendre pairs of length $\ell$ are known if 
\begin{itemize}
    \item $\ell>2$ is a prime (Legendre sequences), see \cite[Section~3]{FGS2001},
    \item  $2\ell+1$ is a prime power (characteristic sequences of Szekeres difference sets), see \cite[Theorem~5.18]{S2017},
    \item $\ell=2^r-1$, $r\ge 2$ ($m$-sequences),
    see for example \cite[Property~5.5]{GG2005},
    \item $\ell=p(p+2)$ where $p$ and $p+2$ are twin primes (twin-prime sequence), see for example \cite{BW2005,OS2010}.
\end{itemize}
Several other constructions for a finite set of $\ell$ are known \cite{BKKT2022,KK2021,TKBG2021}.
The list of integers less than $200$ for
which the question of existence of binary Legendre pairs remains open is
$$\ell=115, 145, 147, 159, 161, 169, 175, 177, 185, 187, 195,$$ see~\cite{BKKT2022}.

Note that for all these $\ell<200$ there exists an Hadamard matrix of order~$2\ell+2$ constructed without using Legendre pairs of order $\ell$. 
The smallest orders for which the existence of an Hadamard matrix is still open are
$$2\ell+2=668, 716, 892, 1132, 1244, 1388, 1436, 1676, 1772, 1916, 1948, 1964,$$ see~\cite{DGK2014}.

A Legendre pair $(A,B)$ is {\em quaternary} if both sequences $A$ and $B$ are quaternary, that is, $A,B\in \{-1,1,-i,i\}^\ell$, where $i=\sqrt{-1}$. Quaternary Legendre pairs exist also for even $\ell$. For example, for $\ell=2$ we take 
$$A=[1,-1]\quad \mbox{and} \quad B=[1,i],$$ that is,
$$PAF(A,1)=-2\quad \mbox{and}\quad PAF(B,1)=0.$$ 

A {\em quaternary} (or {\em complex}) {\em Hadamard matrix} is an $n\times n$-matrix with entries in $\{-1,1,-i,i\}$ which satisfies
$$H\overline{H}^T=nI_n,$$
where $I_n$ denotes the $n\times n$ identity matrix. The order $n$ of a quaternary Hadamard matrix must be even and it is conjectured that there is a quaternary Hadamard matrix for any even order. 
The smallest undecided order is $n=94$, see \cite{D1993,E2019}.

In Section~\ref{LegHad} we construct quaternary Hadamard matrices of order $2\ell+2$ from quaternary Legendre pairs of length $\ell$ (and thus binary Hadamard matrices of order $4\ell+4$).

In Section~\ref{compressed}, for a prime $p>2$ we provide a construction of a `seed' pair of length $p$ from which we can derive a quaternary Legendre pair
of length 
$$\ell=2p \quad  \mbox{for}\quad p=3,5,7,13,19,31,41.$$
The smallest prime for which this method fails is $p=11$. However, we also provide a quaternary Legendre pair of length $22$ as well as pairs for the remaining even $\ell\le 24$
in Section~\ref{rem}.
The main idea to reduce the search space is a modified PSD-criterion that certain PSDs (Power Spectral Densities) are integers which are sums of two integer squares. 
We conjecture that there is a quaternary Legendre pair for any even length~$\ell$, which is stronger than the (complex) Hadamard conjecture. The smallest open case is $\ell=28$.

\section{Quaternary Legendre pairs and (quaternary) Hadamard matrices}
\label{LegHad}


First note that quaternary Legendre pairs are (almost) balanced.

\begin{lemma}
  Write 
  $$A=[a_0,a_1,\ldots,a_{\ell-1}], \quad B=[b_0,b_1,\ldots,b_{\ell-1}]$$
  and let $(A,B)$ be a quaternary Legendre pair of length $\ell$.
  Put
$$\alpha=\sum_{j=0}^{\ell-1} a_j \quad \mbox{and}\quad \beta=\sum_{j=0}^{\ell-1} b_j.$$
  Then we have 
  $$|\alpha|^2+|\beta|^2=2,$$
  $$\alpha,\beta\in \{-1,1,-i,i\}\quad \mbox{if $\ell$ is odd}$$
  and 
  $$\{\alpha,\beta\}\in \{ \{0,1+i\}, \{0, 1-i\},\{0,-1+i\}, \{0,-1-i\}\}\quad \mbox{if $\ell$ is even}.$$
\end{lemma}
Proof. 
We have 
$$|\alpha|^2+|\beta|^2=\sum_{s=0}^{\ell-1}(PAF(A,s)+PAF(B,s))=2\ell+(\ell-1)(-2)=2.$$
For odd $\ell$ we have $\alpha \beta\ne 0$ and since $\alpha,\beta\in \mathbb{Z}[i]$ we get $|\alpha|=|\beta|=1$ and the result follows.
For even $\ell$ we have either $\alpha=0$ and $|\beta|^2=2$ or $|\alpha|^2=2$ and $\beta=0$ and the result follows. \hfill $\Box$\\

For odd $\ell$ we may assume $\alpha=\beta=1$ since otherwise we multiply the sequence~$A$ element-wise by $\overline{\alpha}$ and $B$ element-wise by $\overline{\beta}$ without changing the PAFs.
For even $\ell$ we may assume 
\begin{equation}\label{balance}
\alpha=0\quad  \mbox{and}\quad \beta=1+i
\end{equation}
since otherwise we may interchange~$A$ and~$B$ if $\beta=0$ to get $\alpha=0$ and $|\beta|^2=2$, then we may multiply~$B$ element-wise by $-1$ if $\beta=-1\pm i$ to get $\beta=1\mp i$, and take the complex conjugates of $A$ and $B$ if $\beta=1-i$.

We use the notation
$$C(A)=\left(\begin{array}{ccccc} a_0 & a_1 & \ldots & a_{\ell-2} & a_{\ell-1}\\
a_{\ell-1} & a_0 & \ldots & a_{\ell-3} & a_{\ell-2}\\
\vdots & \vdots &  & \vdots & \vdots\\
a_2 & a_3 & \ldots &a_0 & a_1\\
a_1 & a_2 & \ldots & a_{\ell-1}& a_0\end{array}\right),\quad A=[a_0,a_1,\ldots,a_{\ell-1}],$$
and recall that $C(A)C(B)=C(B)C(A)$, see for example \cite[Theorem~3.2.4]{D1979}. Then the following construction of quaternary Hadamard matrices from quaternary Legendre pairs can be easily verified.

\begin{theorem}
Let $(A,B)$ be a quaternary Legendre pair of odd length $\ell$ with
$$\alpha=\beta=1.$$
Then 
$$\left(\begin{array}{cc|cccccc}
-1 & -1 & 1 & \ldots & 1 & 1 &\ldots & 1\\
-1 & 1 & 1 & \ldots & 1 & -1 &\ldots &-1\\\hline
1 & 1 & & & & & &\\
\vdots & \vdots & & C(A) & & &C(B)&\\
1 & 1 & & & & & &\\\hline
1 & -1 & & & & & &\\
\vdots & \vdots & & C(\overline{B})^T & & &-C(\overline{A})^T&\\
1 & -1 & & & & & &
\end{array} \right)$$
is a quaternary Hadamard matrix of order $2\ell+2$.

Let $(A,B)$ be a quaternary Legendre pair of even length $\ell$ with 
$$\alpha=0\quad \mbox{and}\quad  \beta=1+i.$$
Then 
$$\left(\begin{array}{cc|cccccc}
-1 & i & 1 & \ldots & 1 & 1 &\ldots & 1\\
-i & 1 & 1 & \ldots & 1 & -1 &\ldots &-1\\\hline
1 & 1 & & & & & &\\
\vdots & \vdots & & C(A) & & &C(B)&\\
1 & 1 & & & & & &\\\hline
1 & -1 & & & & & &\\
\vdots & \vdots & & C(\overline{B})^T & & &-C(\overline{A})^T&\\
1 & -1 & & & & & &
\end{array} \right)$$
is a quaternary Hadamard matrix of order $2\ell+2$.
\end{theorem}

Note that a binary Hadamard matrix of order $4\ell+4$ can be constructed from a quaternary Hadamard matrix of order $2\ell+2$,
see \cite[Theorem 1]{C1965}, that is, from a quaternary Legendre pair of length $\ell$.

\section{$2$-compressed pairs and `seeds'}
\label{compressed}

For binary sequences the concept of compression was introduced in \cite{DK2015}. We adjust it here to quaternary sequences of length $\ell=2p$.
Write
$$A=[a_0,a_1,\ldots,a_{2p-1}],\quad B=[b_0,b_1,\ldots,b_{2p-1}]$$
and let $(A,B)\in \{-1,1,-i,i\}^{4p}$ 
be a quaternary Legendre pair of length $2p$. 

The {\em $2$-compression of $(A,B)$} is the pair
$(A_p,B_p)$ of sequences of length $p$ with 
$$A_p=\left[a_0^{(p)},a_1^{(p)},\ldots,a_{p-1}^{(p)}\right],\quad B_p=\left[b_0^{(p)},b_1^{(p)},\ldots,b_{p-1}^{(p)}\right]$$
defined by 
$$a_j^{(p)}=a_j+a_{j+p}, b_j^{(p)}=b_j+b_{j+p}\in \{-2,0,2,-2i,2i,-1-i,-1+i,1-i,1+i\}$$ 
for $j=0,1,\ldots,p-1$.
For prime $p>2$, we start with a theoretic construction of a `seed' pair $(A_p,B_p)$ and 'hope' to be able to decompress it to a quaternary Legendre pair $(A,B)$.
In particular, the decompression $A$ of $A_p$ is essentially unique (up to four choices of $(a_0,a_p)$) and can be explicitly described. So we are only left with finding a decompression $B$ of $B_p$.

\begin{definition}\label{Definition:SeedSequences3}
For a prime $p>2$ 
let two `seed' sequences 
$$A_p=\left[a_0^{(p)},a_1^{(p)},\ldots,a_{p-1}^{(p)}\right]\quad \mbox{and}\quad B_p=\left[b_0^{(p)},b_1^{(p)},\ldots,b_{p-1}^{(p)}\right]$$
be defined by
$$a_j^{(p)}=\left\{\begin{array}{ll} 0, & j= 0,\\
2\left(\frac{j}{p}\right), & j=1,2,\ldots,p-1,\end{array}\right.$$
and
$$b_j^{(p)}=\left\{\begin{array}{ll} 1+i, & j= 0,\\
0,& j=1,2,\ldots,p-1,\end{array}\right.$$
periodically continued with period $p$, 
where $\left(\frac{.}{.}\right)$ is the Legendre symbol.
\end{definition}
Note that the definition of $A_p$ is inspired by the construction of binary Legendre pairs of length $p$ which can be also considered the core of Paley's construction of Hadamard matrices.\\

The {\em discrete Fourier transform} $DFT(A,s)$ of $A$ at lag $s$ is
$$DFT(A,s)=\sum_{j=0}^{\ell-1} a_j \xi_\ell^{js},\quad s=0,1,2,\ldots,\ell-1,$$
where $\xi_\ell$ is the primitive $\ell$th root of unity 
$$\xi_\ell=\exp\left(\frac{2\pi i}{\ell}\right).$$
The {\em power spectral density} $PSD(A,s)$ of $A$ at $s$ is
$$PSD(A,s)=\left|DFT(A,s)\right|^2
=\sum_{j=0}^{\ell-1}PAF(A,j)\xi_\ell^{-js}=\sum_{j=0}^{\ell-1}\overline{PAF(A,j)}\xi_\ell^{js}.$$
For a Legendre pair $(A,B)$ of length $\ell$ we immediately get
$$PSD(A,s)+PSD(B,s)=2\ell+2,\quad s=1,2,\ldots,\left\lfloor \frac{\ell}{2}\right\rfloor.$$
The {\em PSD-criterion} to search for (binary) Legendre pairs $(A,B)$ was introduced in \cite{FGS2001} and
discards all candidate sequences $B$ with 
\begin{equation}
PSD(B,s)>2\ell+2\quad \mbox{for some }s\in\left\{1,2,\ldots,\left\lfloor \frac{\ell}{2}\right\rfloor\right\}
\label{PSDtest}
\end{equation} 
from the search.
Here we can use a much stronger PSD-criterion since we know all values $PSD(B,s)$, $s=1,2,\ldots,\left\lfloor \frac{\ell}{2}\right\rfloor$. 

In the following theorem we determine the four decompressions $A$ of $A_p$ and mention many properties of possible decompressions~$B_p$ of $B$ including all of their PSDs.

\begin{theorem}\label{Thm2}
The construction of Definition~\ref{Definition:SeedSequences3}
satisfies:
\begin{enumerate}
\item $$PAF(A_p,0)=4(p-1),\quad  PAF(B_p,0)=2,$$
 $$PAF(A_p,s)=-4, \quad PAF(B_p,s)=0, \quad s=1,2,\ldots,p-1.$$
\item $$DFT(A_p,s)=\left\{\begin{array}{cc} 2\left(\frac{s}{p}\right)p^{1/2}, & p\equiv 1\bmod 4,\\
2i\left(\frac{s}{p}\right)p^{1/2}, &p\equiv 3\bmod 4,\end{array}\right.$$
$$DFT(B_p,s)=1+i,\quad s=1,2,\ldots,p-1,$$
$$PSD(A_p,s)=4p,\quad PSD(B_p,s)=2,\quad s=1,2,\ldots,p-1,$$
$$PSD(A_p,s)+PSD(B_p,s)=2(2p+1),\quad s=1,2,\ldots,p-1.$$
\item There are exactly four $2$-decompressions $A=[a_0,a_1,\ldots,a_{2p-1}]$ of $A_p$ defined by
$$a_p=-a_0,\quad a_0\in \{-1,1,-i,i\}$$
and
$$a_j=a_{j+p}=\left(\frac{j}{p}\right),\quad j=1,2,\ldots,p-1.$$
They satisfy
$$PAF(A,0)=2p,\quad PAF(A,p)=2p-4,$$
$$PAF(A,s)=PAF(A,2p-s)=-2,\quad s=1,2,\ldots,p-1,$$
$$DFT(A,2s)=2\left(\frac{s}{p}\right)\left\{\begin{array}{cc} p^{1/2}, & p\equiv 1\bmod 4,\\
ip^{1/2}, &p\equiv 3\bmod 4,\end{array}\right.$$
$$PSD(A,2s)=4p,\quad s=1,2,\ldots,p-1,$$
$$DFT(A,2s-1)=2a_0,\quad PSD(A,2s-1)=4,\quad s=1,2,\ldots,p.$$
\item Assume that $B$ is a $2$-decompression of $B_p$. Then we have
$$b_0+b_p=1+i\quad \mbox{and}\quad b_{j+p}=-b_j,\quad j=1,2,\ldots,p-1,$$
$$PAF(B,0)=2p,\quad PAF(B,p)=-2p+2,$$
$$DFT(B,2s)=1+i,\quad PSD(B,2s)=2,\quad s=1,2,\ldots,p-1.$$
\item In addition, if $(A,B)$ is a Legendre pair, we have
$$PAF(B,s)=PAF(B,2p-s)=0,\quad s=1,2,\ldots,p-1,$$
$$PSD(B,2s-1)=4p-2,\quad s=1,2,\ldots,p,$$
$$DFT(B,p)=a+ib \quad\mbox{with integers $a,b$ satisfying }a^2+b^2=4p-2$$
and
$$a\equiv b\equiv 1\bmod 2.$$
\end{enumerate}
\end{theorem}
Proof. 
1. We have 
$$PAF(A_p,s)=4\sum_{j=0}^{p-1}\left(\frac{j}{p}\right)\left(\frac{j+s}{p}\right)
=\left\{\begin{array}{cc} 4(p-1), & s=0,\\
-4, & s=1,\ldots,p-1,\end{array}\right.$$
see for example \cite[Lemma~6.4.7]{NW2015}. 
The values of $PSD(B_p,s)$ are trivial.

2. Next we have
$$DFT(A,2s)=DFT(A_p,s)=2\sum_{j=1}^{p-1}\left(\frac{j}{p}\right)\xi_p^{js}
=2G_s,$$
where $G_s$ is a Gauss sum, that is,
$$G_s=\left\{\begin{array}{cc}\left(\frac{s}{p}\right) p^{1/2}, &p\equiv 1\bmod 4,\\
i\left(\frac{s}{p}\right) p^{1/2}, &p\equiv 3\bmod 4,\end{array}\right.$$
by \cite[Theorems 5.12 and 5.15]{LN1997}. 
The rest of 2.\ is trivial.

3. We only verify
$$DFT(A,2s-1)=\sum_{j=0}^{p-1}(a_j-a_{j+p}) \xi_{2p}^{(2s-1)j}=2a_0$$
and the rest is straightforward.

4. is trivial.

5. Now assume that $(A,B)$ is a Legendre pair, that is in particular,
$$PSD(A,p)+PSD(B,p)=4p+2$$
and thus $PSD(B,p)=4p-2$.
Since 
$$DFT(B,p)\in \mathbb{Z}[i]$$
and
$$DFT(B,p)=\sum_{j=0}^{p-1}(b_j-b_{j+p})
= b_0-b_p+2\sum_{j=1}^{p-1}b_j
\equiv 1+i \bmod 2,$$
that is $DFT(B,p)=a+ib$ for some integers $a\equiv b\equiv 1\bmod 2$,
we have $4p-2=PSD(B,p)=|DFT(B,p)|^2=a^2+b^2$
and the rest is trivial. \hfill $\Box$\\

Remarks. 
1. If we assume that besides $b_{j+p}=-b_j$, $j=1,2,\ldots,p-1$, we also have
\begin{equation}\label{asym}b_{p-j}=-b_j, \quad j=1,2,\ldots,\frac{p-1}{2},
\end{equation}
then we get 
$$b_{2p-j}=b_j,\quad j=1,2,\ldots,p-1,$$
and thus
$$DFT(B,s)=b_0+(-1)^sb_p+\sum_{j=1}^{p-1}b_j\left(\xi_{2p}^{js}+\xi_{2p}^{-js}\right),\quad s=1,2,\ldots,p-1,$$
that is,
$$DFT(B,s)=DFT(B,2p-s),\quad s=1,2,\ldots,p-1,$$
and
\begin{equation}\label{PSDsym}PSD(B,s)=PSD(B,2p-s),\quad s=1,2,\ldots,p-1.
\end{equation}
(Note that without this condition we only have
$$DFT(B,s)=\overline{DFT(\overline{B},2p-s)},\quad s=1,2,\ldots,p-1,$$
and
$$PSD(B,s)=PSD(\overline{B},2p-s),\quad s=1,2,\ldots,p-1.)$$
This has been a suitable condition to reduce the search space.
In addition we may choose 
\begin{equation}\label{b0}
b_0=1\quad\mbox{and} \quad b_p=i.
\end{equation}
In this case we have
$$DFT(B,p)=1-i+4\sum_{j=1}^{(p-1)/2}(-1)^jb_j=a+ib\quad \mbox{with }a^2+b^2=4p-2,$$
and thus $a\equiv -b\equiv 1\bmod 4$,
which further reduces the number of possible $b_1,b_2,\ldots,b_{(p-1)/2}\in \{-1,1,-i,i\}$.

2. Note that $4p-2$ cannot be a sum of two integer squares if $2p-1$ is divisible by an odd power of a prime equivalent to $3$ modulo $4$ by the sum of two squares theorem.
Hence, our first method fails for $p=11,17,29,47,\ldots$

However, it was successful for $p=3,5,7,13,19,31,41$.\\

The decompressed sequences $B=[b_0,b_1,\ldots,b_{2p-1}]$ which lead to quaternary Legendre pairs $(A,B)$ are defined by 
$$b_0=1, b_p=i, \quad b_{j+p}=b_{p-j}=-b_j, \quad j=1,2,\ldots,p-1,$$
and $[b_1,b_2,\ldots,b_{(p-1)/2}]$  given in the following table:
$$\begin{array}{c|c}
p & [b_1,b_2,\ldots,b_{(p-1)/2}]\\\hline
3 & [1]\\
5& [1,i]\\
7&  [ 1, -i, -1] \\
13&  [ 1, 1, -1, i, 1, i]\\
19& [ 1, -1, 1, -i, -1, -i, -i, 1, 1]\\
31& [1, -1, 1, -1, -1, -i, -1, -i, -i, -1, i, i, -i, -1,
-1] \\
41& [1, 1, -i, 1, -1, -i, -1, i, 1, i, 1, -i, i, i, i, -i, 1, 1, -1, -i]
\end{array}$$

\section{The remaining even $\ell\le 24$}
\label{rem}

Again we may assume \eqref{balance}
and denote by 
$$A_k=\left[a_0^{(k)},a_1^{(k)},\ldots,a_{k-1}^{(k)}\right]$$
the {\em $\frac{\ell}{k}$-compression} of $A$, that is,
$$a_j^{(k)}=\sum_{n=0}^{\ell/k-1}a_{kn+j},\quad j=0,1,\ldots,k-1.$$
The alphabet of the $\frac{\ell}{k}$-compressions $A_k$ of $A\in \{-1,1,-i,i\}^\ell$ has $\left(\frac{\ell}{k}+1\right)^2$ possible elements, namely those $a+ib$ with 
$$|a|+|b|\le \frac{\ell}{k}\quad \mbox{and}\quad  a+b\equiv \frac{\ell}{k}\bmod 2,$$ that is, $a_j^{(k)}$, $j=0,1,\ldots,k-1$, is of the form
$$a+i\left(\frac{\ell}{k}-a-2j\right),\quad a=-j,-j+1,\ldots,\frac{\ell}{k}-j, \quad j=0,1,\ldots,\frac{\ell}{k}.$$
Recall
$$PSD\left(A,\frac{s\ell}{k}\right)=PSD(A_k,s),\quad s=1,2,\ldots,k-1.$$
We state now several restrictions on the PSDs which reduce our search space for Legendre pairs.

\begin{theorem}\label{comprcond}
Let $(A,B)$ be a quaternary Legendre pair of length $\ell$ with \eqref{balance}.\\
1. For $\ell\equiv 0\bmod 2$, $PSD(B_2,1)$ is an integer with 
$$PSD(B_2,1)\equiv 2\bmod 8$$ and the square-free parts of both 
$$PSD(B_2,1)\quad \mbox{and}\quad PSD(A_2,1)=2\ell+2-PSD(B_2,1)$$
are not divisible by a prime $\equiv 3\bmod 4.$\\
2. For $\ell\equiv 0\bmod 4$, $PSD(B_4,1)$ is an integer with 
$$PSD(B_4,1)\equiv 2\bmod 8$$ 
and the square-free parts of both
$$PSD(B_4,1)\quad\mbox{and}\quad PSD(A_4,1)=2\ell+2-PSD(B_4,1)$$ 
are not divisible by a prime $\equiv 3\bmod 4.$
\end{theorem}
Proof.
1. We have
$$DFT(B,\ell/2)=\sum_{j=0}^{\ell-1}b_j (-1)^j=a+ib\in \mathbb{Z}[i]$$
for some integers $a,b$.
Since 
$$1+i=\sum_{j=0}^{\ell-1}b_j\equiv a+ib\bmod 2$$ 
we must have
$$a\equiv b\equiv 1\bmod 2$$
and thus
$$PSD(B_2,1)=PSD(B,\ell/2)=a^2+b^2\equiv 2\bmod 8.$$
Similarly, we get
$$PSD(A_2,1)=PSD(A,\ell/2)=c^2+d^2,\quad c,d\in \mathbb{Z}.$$
$PSD(A_2,1)+PSD(B_2,1)=2\ell+2$
and the sum of two squares theorem finishes the proof of the first part.\\

2. First note that $DFT(A,\ell/4),DFT(B,\ell/4)\in \mathbb{Z}[i]$.
Then we have
$$DFT(B,\ell/4)=\sum_{j=0}^{\ell-1}b_ji^j\equiv \sum_{j=0}^{\ell-1}b_j(-1)^j\equiv DFT(B,\ell/2)\equiv 1+i \bmod (1+i)$$
since $DFT(B,\ell/2)\equiv 1+i\bmod 2$ and $1+i$ divides $2=(1+i)(1-i)$ in $\mathbb{Z}[i]$.
Moreover, we get
$$DFT(B,\ell/4)\equiv \sum_{j=0}^{\ell-1}b_j\equiv 1+i \bmod (1-i)$$
and thus
$$DFT(B,\ell/4)\equiv 1+i\bmod 2$$
since $1+i$ and $1-i$ are primes in $\mathbb{Z}[i]$.
Now the second result follows analogously to the first one.
\hfill $\Box$\\

The above condition on $PSD(A,\ell/2)$ and $PSD(B,\ell/2)$ reduces our search to sequences with the following values of these PSDs:

$$\begin{array}{c|c}
\ell & (PSD(A,\ell/2),PSD(B,\ell/2))\\\hline
2 & {\bf (4,2)}\\
4 & {\bf (0,10)},(8,2)\\
8 & (0,18),(8,10),{\bf (16,2)}\\
12 &  (0,26), (8,18), {\bf (16,10)}\\
16& (0,34), (8,26), (16,18), {\bf (32,2)}\\
18 & (4,34), {\bf (20,18)}, (36,2)\\
20 & (8,34), {\bf (16,26)}, (32,10), (40,2)\\
22 & (20,26),{\bf (36,10)}\\
24 & {\bf (0,50)}, (16,34), (32,18), (40,10)\\
\end{array}
$$
The boldfaced pairs correspond to the Legendre pairs $(A,B)$ given in the following table:

$$\begin{array}{c|c}
\ell & A\\\hline
2 & [1,-1]\\
4 & [1,1,-1,-1]\\
8 & [1, 1, -1, -i, -1, 1, -1, i] \\
12 & [1, 1, 1, -1, i, -i, 1, -1, -1, -1, -i, i] \\
16 &   
 [1, 1, 1, -1, -i, -1, i, 1, i, -1, i, -1, -i, -i, i, -i] 
\\
18 & [1, 1, -1, i, -1, -1, -i, 1, -1, i, -i, 1, -1, -i, 1, -i, i, i] \\
20& [1, 1, 1, -1, i, -i, -i, i, i, 1, i, -1, -1, -i, -1, 1, -i, -i, i, -1]
 \\
22 & [1, 1, 1, 1, -1, -1, -1, 1, -1, 1, -1, -i, i, 1, -1, -1, -i, i, -i, -i, i, i]\\ 
24 & [1, 1, 1, 1, 1, -1, 1, 1, -1, -1, 1, -1, 1, -1, -1, 1, -1, 1, -1, -1, -1, 1, -1, -1]
\end{array}$$
$$\begin{array}{c|c}
\ell &  B\\\hline
2 & [1,i]\\
4 & [1,i,1,-1]\\
8 &   [1, 1, 1, i, -1, 1, -1, -1]\\
12 &  
[1, 1, -1, -1, 1, -1, 1, -1, -1, 1, 1, i]\\
16 &   
[1, 1, 1, -1, 1, 1, -1, 1, -1, i, -i, -i, -1, i, i, -1]
\\
18 &  [1, 1, -1, i, -i, -1, 1, -1, 1, i, -1, -1, i, i, 1, 1, -i, -i]\\
20& 
[1, 1, 1, -i, -i, -i, i, -1, -i, i, -i, i, 1, -1, i, 1, i, -1, i, -1]\\
22 &[1, 1, 1, 1, -1, 1, 1, -1, -i, 1, -i, i, i, i, -1, -i, -1, i, -i, -1, i, -1]\\
24 & 
[1, 1, 1, 1, 1, i, -1, -1, -i, 1, -1, i, 1, -1, -i, -1, 1, i, -1, 1, -i, -1, -1, i]
\end{array}$$

For $\ell\equiv 0\bmod 4$ the set of possible pairs of values $(PSD(A,\ell/4),PSD(B,\ell/4))$ is the same as for $(PSD(A,\ell/2),PSD(B,\ell/2))$ but the values can be different. Here is a list of these values for the Legendre pairs in the above table:
$$\begin{array}{c|cccccc}
\ell & 4 & 8 & 12 & 16 & 20 & 24  \\\hline
(PSD(A,\ell/4),PSD(B,\ell/4)) & (8,2) & (8,10) & (16,10) & (16,18) & (32,10) & (0,50) 
\end{array}
$$

Remarks.
1. For $\ell\equiv 0\bmod 6$, if we assume that the $\frac{\ell}{3}$-compression $A_3$ 
satisfies
$$a_0^{(3)}-a_2^{(3)}, a_1^{(3)}-a_2^{(3)}\in \mathbb{Z},$$
then 
$$PSD(A_3,1) \in \mathbb{Z},$$
since $DFT(A_3,1)=(a_0^{(3)}-a_2^{(3)})+(a_1^{(3)}-a_2^{(3)})\xi_3$.
In this case $PSD(A,\ell/3)=PSD(A_3,1)=PSD(A_6,2)$ is of the form $(a+b\xi_3)(a+b\xi_3^2)=a^2-ab+b^2$, $a,b\in \mathbb{Z}$, and thus its square-free part is not divisible by a prime $\equiv 2\bmod 3$, see for example \cite[Chapter~9.1]{IR1990}.
Moreover, if $A_6$ satisfies
$$a_0^{(6)}-a_3^{(6)}+a_5^{(6)}-a_2^{(6)},a_4^{(6)}-a_1^{(6)}+a_5^{(6)}-a_2^{(6)}\in \mathbb{Z},$$
then the square-free part of $PSD(A,\ell/6)=PSD(A_6,1)$ is not divisible by a prime $\equiv 2\bmod 3$.

Note that both conditions on $A_3$ and $A_6$ are fulfilled if $A_6$ is integral.
For example, our constructions of $A$ for $\ell=6$, that is, $A=A_6=[1,1,-1,-1,1,-1]$, for $\ell=12$, that is, 
$A_6=[2,0,0,-2,0,0]$ and $\ell=24$, that is, 
$A_6=[2,0,-2,2,0,-2]$, are all integral and thus satisfy these conditions.
In particular, we have 
$$(PSD(A_6,1),PSD(A_6,2))=(4,12)\mbox{ for }\ell=6,$$
$$PSD(A_6,1)=PSD(A_6,2)=16\mbox{ for }\ell=12,$$
and
$$(PSD(A_6,1),PSD(A_6,2))=(0,48)\mbox{ for }\ell=24.$$


2. If $a_2^{(3)}=a_1^{(3)}$, then $DFT(A_3,1)\in \mathbb{Z}[i]$ and $PSD(A_3,1)=PSD(A_6,2)$ is a sum of two squares, that is, its square-free part is not divisible by a prime $p\equiv 3\bmod 4$.
If $a_4^{(6)}+a_5^{(6)}=a_1^{(6)}+a_2^{(6)}$, then also $PSD(A_6,1)$ is not divisible by a prime $p\equiv 3\bmod 4$.

Both conditions are fulfilled for our solution for $\ell=12$.
The second condition is also fulfilled for our solutions for $\ell=6$ and $\ell=24$.

3. For $\ell\equiv 0\bmod 6$, there are further very promising `seeds' of the form $A_3=[0,a+ib,-(a+ib)]$
with integers $a$ and $b$ satisfying $|a|+|b|\le \ell/3$, $a^2+b^2\le (2\ell+2)/3$ and $|a|+|b|\equiv \ell/3 \bmod 2$.
In this case we have
$$DFT(A,\ell/3)=DFT(A_3,1)=(a+ib)(\xi_3-\xi_3^2)=(a+ib)\sqrt{-3}$$
and thus
$$PSD(A_3,1)=3(a^2+b^2).$$
For example, our quaternary Legendre pair for $\ell=18$ satisfies 
$$A_3=[0,1-i,-1+i]\quad \mbox{and} \quad PSD(A_3,1)=6.$$


\section{Computational and programming details}

In order to conduct a systematic investigation of 2-decompression for quaternary Legendre pairs using the pairs of `seed' sequences of Definition~\ref{Definition:SeedSequences3}, we wrote a Maple meta-program that produces efficient C code that implements the 2-decompression for `seed' sequences with elements from $\{-2,0,+2,1+i,-1-i\}$. Using this meta-program we found the examples in Section~\ref{compressed}. To give some idea of the complexity involved in these computations, we mention the relevant information for the pattern of Definition~\ref{Definition:SeedSequences3}, namely:
\begin{equation}
    \label{pattern_defn_3}
    A_p = \left[0,2\left(\frac{1}{p}\right),2\left(\frac{2}{p}\right),\ldots,2\left(\frac{p-1}{p}\right)\right], \quad B_p = [1+i, \underbrace{0 \ldots 0}_{p-1 \mbox{ terms}}].
\end{equation}
There are exactly four $2$-decompressions of $A_p$
and $2 \cdot 4^{p-1} = 2^{2p-1}$ possible $2$-decompressions of $B_p$, see Theorem~\ref{Thm2}. All four $2$-decompressions of $A_p$ give the same PAFs and PSDs. 
We restricted ourselves to the cases \eqref{asym} and \eqref{b0}
to reduce the search space for decompressions to $4^{(p-1)/2}=2^{p-1}$ possible choices of $b_1,b_2,\ldots,b_{(p-1)/2}\in \{-1,1,-i,i\}$.

Now by \eqref{PSDsym} and since $PSD(B,2s)=2$, $s=1,2,\ldots,p-1$, for any decompression $B$ of $B_p$, for any $(b_1,b_2,\ldots,b_{(p-1)/2})\in \{-1,1,-i,i\}^{(p-1)/2}$
we have to check the condition 
$$PSD(B,2s-1)=4p-2,\quad s=1,2,\ldots,\frac{p+1}{2},$$
where
$$PSD(B,2s-1)=|DFT(B,2s-1)|^2=\left|1-i+4\sum_{j=1}^{(p-1)/2}b_j\cos\left(\frac{ (2s-1)j\pi}{p}\right)\right|^2.$$
(Note that this procedure can only succeed if the square-free part of $2p-1$ is not divisible by a prime $\equiv 3\bmod 4$ since otherwise $PSD(B,p)=4p-2$ is not possible.)

%

For the sequences found in Section~\ref{rem}, instead of using the classical PSD-test (\ref{PSDtest}), we use a much more precisely targeted PSD test, given that we know explicitly all the PSD values of the desirable decompressed sequence. 


First we determined all eligible pairs 
$$(PSD(A_2,1),PSD(B_2,1)).$$ For $\ell\equiv 0\bmod 4$ the same values are eligible for $(PSD(A_4,1),PSD(B_4,1))$. Then we determined all compressions satisfying these conditions and decompressed them.
Some additional ideas mentioned in the remarks also helped to reduce the search space.

Finally, note that in contrast to the case $\ell=24$ the case $\ell=26$ was settled with the more efficient method of Section~\ref{compressed}. Compared to $\ell=24$, it seems to be natural to assume that the complexity of finding a Legendre pair for $\ell=28$ increases by the factor $256$. This may be in reach with some more programming effort such as parallelization or using faster computers. However, we stopped here and prefer to search for suitable `seeds' which may cover more cases, say, of the form $\ell=4p$.

\end{document}